\documentclass[12pt]{article}
\usepackage[cp1250]{inputenc}
\usepackage{amsmath}
\usepackage{amssymb}
\usepackage{amsmath}
\usepackage{amsthm}
\usepackage{epsfig}
\usepackage{graphics}
\usepackage{color}
\usepackage{cancel}
\usepackage{enumerate}
\usepackage{psfrag}
\input{trans}

\def\RR{\mathbb R}

\def\S3{{$\mathbb S^3$}}

\def\area{\operatorname{area}}

\newcommand\const{\operatorname{const}}
\newcommand\enpr{\vbox to 7pt{\hbox to 7pt{\vrule height 7pt width 7pt}}}

\newtheorem{ex}{Example}

\newtheorem{prop}{Proposition}
\newtheorem{thm}{Theorem}

\title{Dupin cyclides osculating surfaces}
\author{Adam Bartoszek, Pawe\l\/ G. Walczak, Szymon M. Walczak {\thanks{The second and the third author were supported by the Polish NSC grant N° 6065/B/H03/2011/40}}\\Wydzia\l\/ Matematyki i
Informatyki\\ Uniwersytet \L{\'o}dzki\\ \L{\'o}d\'z, Poland\\ E-mail:
mak, pawelwal, sajmonw@math.uni.lodz.pl} 

\begin{document}
\maketitle
\section{Introduction}
According to {\it Longman Dictionary of Contemporary English}, "{\it osculation}" means "the act of kissing". In mathematics however, 
"{\it to osculate}" means rather "to be tangent" each to the other, in fact tangent "as much as possible", that is to possess the highest possible degree of tangency. So, generically, a space curve has its osculating plane and osculating sphere, a surface in the 3-dimensional space has -- generically again -- two osculating spheres and so on. Also, osculation of other standard surfaces, for example, quadrics, to general surfaces was of some interest at least since $19^{th}$ century (see \cite{da1}, \cite{dr1}, \cite{dr2} etc.). Let us note that the interest in problems of this type grew recently due to applications in computer aided geometric design (CADG).

The notion of {\it osculating sphere} belongs to extrinsic conformal geometry: if you transform a surface by a conformal (that is, M\"obius) transformation of the 3-dimensional space $\RR^3$(or, better, of the 3-dimensional sphere \S3 seen as the 1-point compactification of $\RR^3$), then this transformation maps the spheres osculating the original surface onto the spheres osculating its image: in fact, M\"obius transformations of \S3 (resp., of $\RR^3$) map spheres onto spheres (resp., planes or spheres onto planes or spheres) and -- being diffeomorphisms -- preserve order of tangency.

\begin{figure}
\begin{center}
\includegraphics[scale=.5]{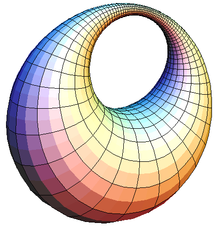}
\caption{A Dupin cyclide\label{fig:cyclide}}
\end{center}
\end{figure}

Another conformally invariant class of surfaces consists of {\it canal surfaces}, that is envelopes of one-parameter families of spheres. These are characterized locally by vanishing of one of their {\it conformal principal curvatures} (see Section 
\ref{sec:lci}). Canal surfaces are interesting from several pints of view; they play some role in computer graphics, medicine etc. Among them, one can distinguish {\it special canals} (see \cite{blw}) characterized by the following property: one of their conformal principal curvatures vanishes identically while the other one is constant along their {\it characteristic circles}, that is the lines of curvature corresponding to the first principal curvature. A smaller class -- again, conformally invariant -- of canal surfaces consists of {\it Dupin cyclides} (shortly, {\it cyclides}, see Figure \ref{fig:cyclide}) which are canal surfaces in two ways: they are envelopes of two different one-parameter families of spheres, therefore can be characterized locally by vanishing of both their conformal principal curvatures. Dupin cyclides where defined first in \cite{du} and  studied later on in, for example, \cite{ma} and \cite{ca}. Recently, Dupin cyclides have been revived because of their applications in computer aided geometric design (CADG). 

This article is devoted to the study of cyclides osculating general surfaces. We show (Theorem \ref{thm1}) that generically, at any point of a surface, one  has a one-parameter family of cyclides tangent to a surface curve of order three and among them just one is tangent to this curve of order four. This one will be called the {\it osculating cyclide} here. Directions of this tangency of higher order form a line filed on the surface and its integral curves will be called {\it Dupin lines} on the surface under consideration. (Following Drach, \cite{dr1} and \cite{dr2}, terminology we should speak also about {\it lines of cyclidal osculation} or so.) Our Dupin lines are analogous to classical lines of curvature corresponding in the same to the eigenvectors of the Weingarten (shape) operator of a surface, that is to the directions of tangency of order two of its osculating spheres. Each point of a generic surface meets two orthogonal lines of curvature but only one Dupin line. Our main result (Theorem \ref{thm2}) shows that a large class of foliations of open planar domains can be realized as {\it Dupin foliations}, that is foliations by Dupin lines, on several surfaces.  

\section{Local conformal invariants}\label{sec:lci}

Let $S$ be an oriented  surface in \S3 (or, $\RR^3$).
Assume that $S$ is {\it
umbilic free}, that is, that the principal curvatures $k_1(x)$ and $k_2(x)$
of $S$ are different at any point $x$ of $S$. Let $X_1$ and $X_2$ be unit
vector fields tangent to the curvature lines corresponding to,
respectively, $k_1$ and $k_2$. Throughout the paper, we assume that $k_1
> k_2$. Put $\mu = (k_1 - k_2)/2$. Since more than 100 years, it is known
(\cite{tr}, see also \cite{csw} ) that the vector fields $\xi _i = X_i/\mu$
and the coefficients $\theta _i$ ($i = 1,2$) in
\begin{equation*}
[\xi _1, \xi_ 2] = -\frac{1}{2}\left(\theta _2\xi _1 + \theta _1\xi
_2\right)
\end{equation*}
are invariant under arbitrary (orientation preserving) conformal
transformation of $\mathbb R^3$. (In fact, they are invariant under
arbitrary conformal change of the Riemannian metric on the ambient space.)
Elementary calculation involving Codazzi equations shows that
\begin{equation*}
\theta _1 = \frac{1}{\mu ^2}\cdot X_1(k_1)\quad\text{and}\quad\theta _2 =
\frac{1}{\mu ^2}\cdot X_2(k_2).
\end{equation*}
The quantities $\theta _i$ ($i = 1,2$) are called {\it conformal principal
curvatures} of $S$.

Another conformally invariant scalar quantity $\Psi$ can be derived from
the derivation of Bryant's (see \cite{br}) {\it conformal Gauss map}
$\beta$:
\begin{eqnarray*}
\frac{1}{2}\big( \langle\xi _1(\xi _1(\beta)), \xi _1(\xi
_1(\beta))\rangle - \langle\xi _2(\xi _2(\beta)), \xi _2(\xi
_2(\beta))\rangle\\ - \langle\xi _1(\xi _1(\beta)), \xi _2(\beta)\rangle
^2 + \langle\xi _2(\xi _2(\beta)), \xi _1(\beta)\rangle\big)\\=
\Psi - \frac{1}{2}\left(\theta _1^2 - \theta _2^2 + \xi _1(\theta _1) + \xi
_2 (\theta _2)\right).
\end{eqnarray*}
Note that both sides of the above equality are equal to
\begin{equation*}
\frac{1}{\mu ^3}\left( \triangle H + 2\mu ^2H\right) ,
\end{equation*}
where $H$ is the mean curvature of $S$ and $\triangle$ is the Laplace
operator on $S$ equipped with the Riemannian metric induced from the
ambient space. Moreover, this quantity appears in the Euler-Lagrange
equation for {\it Willmore functional}
\begin{equation*}
\int _S\mu ^2d\area .
\end{equation*}

The vector fields $\xi _1, \xi _2$ (or, the dual 1-forms $\omega _1, \omega
_2$) together with quantities $\theta _1, \theta _2$ and $\Psi$ generate all
the local conformal invariants for surfaces and determine a surface up to
conformal transformations of $\mathbb R^3$ (\cite{fi}, see again
\cite{csw}).

Define $(5\times 5)$ matrices $A_1$ and $A_2$ by
\begin{equation}\label{eq1}
A_1 = \begin{pmatrix}\theta _1/2 & -(1+\Psi)/2 & b/2 & \theta _1/2 & 0\\
1 & 0 & 0 & -1 & (1 + \Psi)/2\\ 0 & 0 & 0 & 0 & - b/2\\
0 & 1 & 0 & 0 & - \theta _1/2\\ 0 & -1& 0 & 0 & - \theta _1/2\end{pmatrix}
\end{equation}
and
\begin{equation}\label{eq2}
A_2 = \begin{pmatrix}-\theta _2/2 & - c/2 & -(1-\Psi)/2 &  \theta _2/2 &
0\\ 0 & 0 & 0 & 0 & c/2 \\ 1 & 0 & 0 & 1 & (1 - \Psi)/2\\ 0 & 0 & -1 & 0 &
- \theta _2/2\\ 0 & 0& -1 & 0 & \theta _2/2\end{pmatrix},
\end{equation}
where $b = -\theta _1\theta _2 + \xi _2(\theta _1)$ and $c = \theta
_1\theta _2 + \xi _1(\theta _2)$.

For further use, let us define also the quantities $a = 3 + \theta_1^2 + \xi_1 (\theta_1 )$ and $d = -3 - \theta_2^2 + \xi_2 (\theta_2)$. 

Given on a simply connected
domain $U\subset\RR^2$ linearly independent 1-forms $\omega _1$ and $\omega
_2$ and smooth functions $\theta _1$, $\theta _2$ and $\Psi$ for which the
matrix valued 1-form $\omega$,
\begin{equation}\label{eq3}
\omega = A_1\omega _1 + A_2\omega _2
\end{equation}
satisfies the structural equation
\begin{equation}\label{eq4}
d\omega + (1/2)[\omega , \omega ] = 0
\end{equation}
there exists an immersion $\iota :U\to\RR^3$ for which $S = \iota (U)$
realizes these forms and functions as local conformal invariants. Let us recall that such $\iota$ is unique up to conformal transformations of the ambient space.

In case of
\begin{align*}
\omega_i=f_idx_i,\ f_i:U\longrightarrow \mathbb R, {f_i}_{\vert U}\not=0\ i=1,2,
\end{align*}
due to symmetries of matrices of the Lie algebra of the M\"obius group, equation (\ref{eq4}) reduces to the following system of four scalar differential equations:

\begin{align}
\partial_2 f_1+(1/2)f_1f_2 \theta_2 &=0,  \label{eq401}\\
\partial_1 f_2-(1/2)f_1f_2\theta_1 &=0,  \label{eq402}\\
f_1 \partial_2 \Psi - f_2 \partial_1 c - f_1f_2 \left[c\theta_1+\theta_2 (\Psi+2 )\right] &=0, \label{eq403}\\
-f_1 \partial_2 b + f_2 \partial_1 \Psi + f_1f_2 \left[ b\theta_2+\theta_1(\Psi-2)\right] &=0. \label{eq404}
\end{align}

Transforming equations (\ref{eq403}) and (\ref{eq404}) with the use of 
$\xi_i = (1/f_i)\cdot (\partial/\partial x_i)$  we obtain
\begin{equation}\xi_2\psi = \xi_1 c+[c\theta_1+\theta_2(\psi +2)],\label{eq4003}
\end{equation}
and
\begin{equation}\xi_1\psi = \xi_2 b-[b\theta_2+\theta_1(\psi -2)]. \label{eq4004}
\end{equation}
Consequently, we have
\begin{align*}
(\xi_1\theta_2&+\xi_2\theta_1)\psi =\\
&-6\theta_1\theta_2 + \frac{3}{2}[b\theta_2^2-c\theta_1^2] - \frac{5}{2}[\theta_2(\xi_2 b)+\theta_1(\xi_1 c)]\\
&+ 2(\xi_2\theta_1-\xi_1\theta_2) + \xi_2^2b-\xi_1^2c-c(\xi_1\theta_1) - b(\xi_2\theta_2).
\end{align*}
and, in a generic situation, that is when
\begin{equation}\label{eq411}
\xi_1(\theta_2) + \xi_2(\theta_1)\ne 0,
\end{equation}
substituting the right hand sides of the formulae expressing $b$ and $c$  (right after formulae (\ref{eq1}) and (\ref{eq2})) we can express $\Psi$ in terms of $\theta_i$'s and their derivatives:  
\begin{align}\label{eq4001}
\psi &= \frac{1}{\xi_1\theta_2+\xi_2\theta_1} [-6\theta_1\theta_2 + 2(\xi_2\theta_1-\xi_1\theta_2)  \nonumber\\
&+ 4(\theta_2^2(\xi_2\theta_1) - \theta_1^2(\xi_1\theta_2)) - \frac{3}{2}(\theta_1\theta_2^3 + \theta_2\theta_1^3) \nonumber\\
&-3(\xi_1\theta_1)(\xi_1\theta_2)-3(\xi_2\theta_1)(\xi_2\theta_2)\nonumber\\
&+\frac{7}{2}\theta_1\theta_2(\xi_2\theta_2-\xi_1\theta_1) -\frac{7}{2}(\theta_2(\xi_2^2\theta_1) + \theta_1(\xi_1^2\theta_2))\nonumber\\
&-\theta_1(\xi_2^2\theta_2)-\theta_2(\xi_1^2\theta_1)
+\xi_2^3\theta_1-\xi_1^3\theta_2]
\end{align}
Note that if condition (\ref{eq411}) is not satisfied then $\Psi$ becomes independent on $\theta_i$'s. For example, 
on Dupin cyclides, one has $\theta_1 = \theta_2 = 0$ while $\Psi$ may be an arbitrary constant.

Note also that (\ref{eq4001}) implies the following: (1) if $\theta_1,\theta_2$ are constant, then
\begin{equation*}
 \theta_1\theta_2(\theta_1^2+\theta_2^2+4)=0
\end{equation*}
as was observed in \cite{bw}; (2) if $\theta_1=0$ (we are on a canal surface), then
\begin{equation*}
\psi = -2 - \frac{\xi_1^3\theta_2}{\xi_1\theta_2}
\end{equation*}
as was observed in \cite{blw}.

\section{Osculating cyclides}
Let us recall after \cite{fi} (see \cite{csw} again) that, given non-umbilical point  $p$ of a surface $S$, $S$  can be mapped by an unique M\"obius transformation $g$ of $\RR^3$ to the {\it standard position}, that is to such a position that $g(p) = (0,0,0)$, the tangent plane $T_{g(p)} g(S)$ coincides with the $(x, y)$ plane and $g(S)$ is given locally by the equation
\begin{eqnarray}\label{eq5}
z &=& \frac{1}{2}(x^2 - y^2) + \frac{1}{6} (\theta_1 x^3 + \theta_2 y^2)\nonumber\\ &+& \frac{1}{24}(ax^4 + 4bx^3y +6\Psi x^2y^2 + 4cxy^3 + dy^4) + O(5),
\end{eqnarray}
$\theta_1$, $\theta_2$ and $\Psi$ being local conformal invariants of $S$ at $p$, $a, b, c, d$ being the quantities (at $p$) defined in Section \ref{sec:lci} and $O(5)$ denoting terms of higher order.

If $S = C$ is a cyclide, (\ref{eq5}) reduces to
\begin{equation}\label{eq6}
z = \frac{1}{2}(x^2 - y^2) + \frac{1}{8}(x^4 - y^4) + \frac{1}{6}\Psi_Cx^2y^2 + O(5),
\end{equation}
$\Psi_C$ being the corresponding conformal invariant (at $p = (0,0,0)$) of $C$.

Consider now a surface $S$ and a cyclide $C$ given locally by (\ref{eq5}) and(\ref{eq6}) along a line $m: y = tx$, $t$ being a constant.
Along $m$, the equations describing $S$ and $C$ reduce, respectively, to
\begin{eqnarray}\label{eq7}
z &=& \frac{1}{2}(1 - t^2)x^2 +\frac{1}{6}(\theta_1 + \theta_2 t^3)x^3\nonumber\\ 
&+& \frac{1}{24}(a + 4bt + 6\Psi t^2 + 4ct^3 + dt^4)x^4 + O(5)
\end{eqnarray}
and
\begin{equation}\label{eq8}
z = \frac{1}{2}(1 - t^2)x^2 + \Big(\frac{1}{8}(1 - t^4) + \frac{1}{6}\Psi_Ct^2\Big)x^4 + O(5).
\end{equation} 

If one of the local conformal curvatures $\theta_i$, say $\theta_2$, is different from zero, then $S$ and $C$ are tangent of order $3$ in the direction of the line $m$ corresponding to the value $t = \sqrt[3]{\theta_1/\theta_2}$. The cyclides tangent of order $3$ at such a point to our surface $S$ are parametrized by the invariant $\Psi_C$ which may take any value in $\RR$. That is, there exists a one-parameter family of cyclides like that. If, moreover, the other conformal principal curvature (here, $\theta_1$) is also different from zero, then among all of them there exists exactly one which is tangent to $S$ at $p$ in the direction of $m$ of order $4$. Let us denote it by ${\cal C}$. The cyclide ${\cal C}$ corresponds to the value of $\Psi_C$ satisfying
\begin{equation}\label{eq9}
\frac{1}{24}(a + 4bt + 6\Psi t^2 + 4ct^3 + dt^4) = \frac{1}{8}(1 - t^4) + \frac{1}{6}\Psi_Ct^2
\end{equation}
with $t = \sqrt[3]{\theta_1/\theta_2}$ and is said to be {\it osculating} here. Observe, that all the cyclides of our one-parameter family have at the reference point $p$ the same principal directions and osculating spheres (which belong to the two one parameter families enveloped by them) as the surface $S$.
 
For further use, we will adopt the following terminology. If at least one of the conformal principal curvatures of a surface $S$ at a point $p$, $p\in S$, equals zero, then $p$ is a {\it canal} (or, {\it ridge}) {\it point} of $S$; if both of them are equal to zero at $p$, then $p$ is a {\it Dupin point}. Points which are not canal are  said to be {\it non-canal}, points which are not Dupin are {\it non-Dupin}. Note that {\it ridges} (that is curves built of ridge points corresponding to one of the equations $\theta_1 = 0$ and $\theta_2 = 0$) are of some interest in the computer graphics
(see, for example, \cite{cg} or \cite{po}).

With this terminology we have the following.

\begin{thm}\label{thm1}
 At any non-Dupin point $p$ of a surface $S$ there exists a one-parameter family of cyclides tangent at $p$ of order $3$ to $S$ in the direction making with the lines of curvature angles $\alpha$ and $\pi/2 - \alpha$ such that $\tan \alpha = \sqrt[3]{\theta_1/\theta_2}$; if $p$ is non-canal point, there exists at $p$ exactly one  osculating cyclide, that is the cyclide tangent to $S$ at $p$ in the direction $\alpha$ of order $4$.
\hfill\qed
\end{thm}

From the above, it follows that if $S$ is a surface without Dupin points, then the vector field
\begin{equation}\label{eq91}
V = \sqrt[3]{\theta_2}\cdot X_1 + \sqrt[3]{\theta_1}\cdot X_2
\end{equation}
is non-singular and defines a 1-dimensional foliation ${\cal O}$. A leaf of this foliation is a curve $\gamma$ on $S$ such that at every its point $\gamma(t)$ the osculating cyclide ${\cal C}(t)$ of $S$ at $\gamma(t)$ is tangent to $\gamma$ of order $4$. The leaves of ${\cal O}$ play the role analogous to this of curvature lines and are called hereafter {\it Dupin lines}; ${\cal O}$ itself is called here the {\it Dupin foliation}.

\begin{figure}
\begin{center}
\includegraphics[scale=.5]{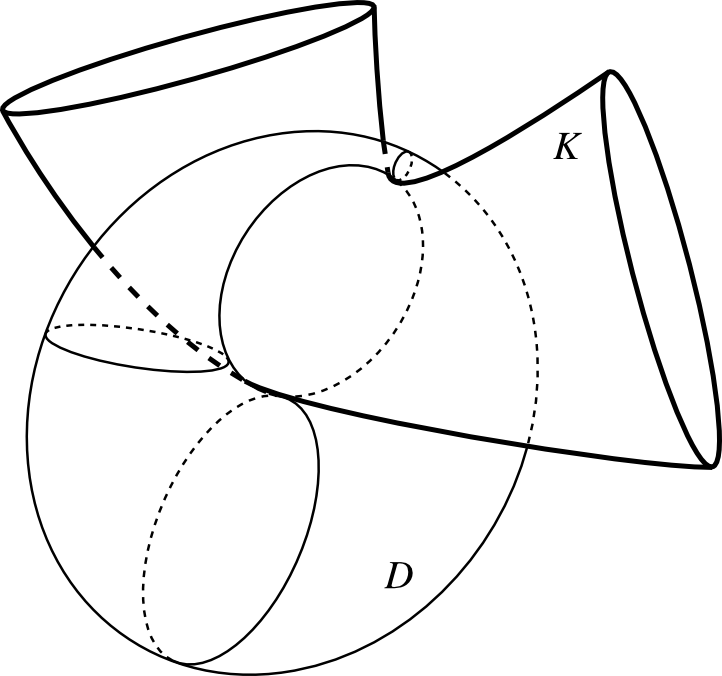}
\caption{Dupin necklace $D$ of a canal $K$\label{Dupin_neck}}
\end{center}
\end{figure}

\begin{ex}\label{ex1} {\rm A (piece of a) canal surface $S$ is characterized by vanishing of one, say $\theta_1$, of its conformal principal curvatures. On such $S$, curvature lines corresponding to the first principal curvature $k_1$ are circles called {\it characteristic circles} of the canal. In this case, Dupin lines coincide with these circles and the osculating cyclides, see Figure \ref{Dupin_neck}, are tangent to $S$ along them. In \cite{blw}, the existence of such cyclides has been established and the name {\it Dupin necklace} has been invented.}
\end{ex} 

\section{Dupin and Darboux lines} 

{\it Darboux lines}{\footnote{We are aware of the fact that the term "Darboux curves" has another meaning in algebraic geometry while our "Darboux lines" are called also "D-curves", see \cite{sa1}, \cite{sa2} etc.; since the letter "D" doesn't distinguish between Darboux and Dupin we decided to use their full names in our terminology.} on a surface $S$ are curves $\gamma$ satisfying the following condition: at any point of $\gamma$, the osculating sphere of $\gamma$ is tangent to $S$. Obviously, the notion of Darboux lines belongs to conformal geometry. Also,
it is known \cite{sa1} that at any non-umbilical point $p$ of $S$ there exists a unique Darboux curve in the direction of a vector $v\in T_pS$ different from the vectors of principal curvatures.  

Recently, R. Garcia, R. Langevin and the third author \cite{glw}, studied the dynamics of the flow determined by Darboux lines. Among the others, they have shown that the angle $\alpha$ between a Darboux line $\gamma$ and lines of curvature satisfies the equation
\begin{equation}\label{eq10}
12\sin\alpha\cos\alpha\frac{d\alpha}{d\sigma} = \theta_1\cos^3\alpha + \theta_2\sin^3\alpha,
\end{equation}
$d\sigma$ being the arc length along $\gamma$ normalized by the factor $(k_1 - k_2)^{-1}$ ($k_1\ne k_2$), so that it becomes conformally invariant.

From the above we can extract directly the following observation.

\begin{prop}\label{prop1}
The leaves of the Dupin foliation ${\cal O}$ on a surface $S$ are tangent to the Darboux lines for which the point of tangency with a leaf is critical for the angle between the Darboux line and the principal directions on $S$. Generically, under the condition $V(\log |\theta_1| + \log |\theta_2|)\ne 0$, $V$ being the vector field defined in (\ref{eq91}), at every such a point $\alpha$ has its local extremum.  \hfill\qed
\end{prop}  

\section{Intersections}\label{sec:caps}
It is known (and not difficult to prove) that -- at a non-umbilical point $p$ of a surface $S$, where the principal curvatures are equal
to $k_1$ and $k_2$ --
the sphere $\Sigma_\alpha$ which is tangent to $S$ at $p$ and has the normal curvature
\begin{equation*}
k = \cos^2\alpha k_1 + \sin^2\alpha k_2, \quad 0 \le\alpha\le \pi/2 ,
\end{equation*}
intersects $S$ along two curves which meet at $p$ at the angle $2\alpha$. In particular, the osculating spheres $\Sigma_0$ and
$\Sigma_{\pi/2}$ intersect $S$ long curves making at $p$ a cusp, while the mean sphere $\sigma_{\pi/4}$ intersects $S$ along two
orthogonal curves. 

Intersections of surfaces $S$ given by (\ref{eq5}) with the Dupin cyclides (\ref{eq6}), in particular of $S$ and the osculating
cyclide,  are also of some interest. Solving the system of equations (\ref{eq5})-(\ref{eq6}) one can observe that the solution has
a number of components, exactly one of them passing through the reference point of $S$. The intersection is simpler (some of these
components disappear) in the case of the osculating cyclide. Figure {\ref{fig21} shows the intersections (more precisely, their 
projections to the common tangent plane) of a generic surface (given by canonical equation (\ref{eq5}) with zero higher order terms)
with the cyclides $C$ (\ref{eq6}) for different values of $\Psi_C$. The thick lines correspond to the osculating cyclide, that is for
the value of $\Psi_C$ which satisfies (\ref{eq9}).

\begin{figure}
\begin{center}
\includegraphics[scale=.3]{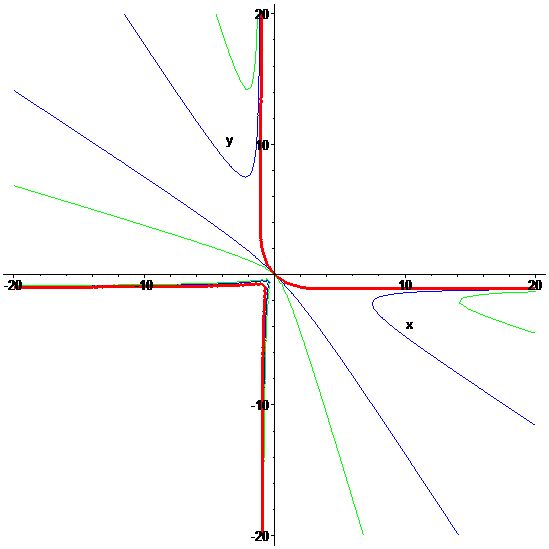}
\caption{Intersections of Dupin cyclides with a surface\label{fig21}}
\end{center}
\end{figure} 

\section{Prescribing Dupin foliations}\label{sec:dupfol}
Assume that we have given a real-valued function $\alpha$ on a simply-connected (say, convex) open domain $U\subset\RR^2$. 
Here, we are looking for a surface $S$ such that the leaves of the Dupin foliation on $S$ intersect the lines of 
curvature at angle $\alpha$. If so, 
the quotient $\theta_1/\theta_2$ should coincide with $\kappa = \tan^3\alpha$.  Assume that $\kappa\ne 0$. 
Equations (\ref{eq401}) and (\ref{eq402}) together with $\theta_1 = \kappa\theta_2$ imply that
\begin{equation}\label{eq101}
 \theta_2 = -2\partial_2f_1/(f_1 f_2) = 2\partial_1f_2/(\kappa f_1 f_2),
 \end{equation}
 therefore that
 \begin{equation}\label{eq102}
 -\partial_2f_1 = \partial_1f_2/\kappa .
 \end{equation}
One can solve (\ref{eq102}) choosing arbitrary (say, positive) $f_2$ and arbitrary values of $f_1$ along a chosen line $x_2 = \const$
and integrating the function $\partial_1f_2/\kappa$ along the lines $x_1 = \const$. Given nonzero (say, positive) functions $f_1$ and $f_2$ on $U$ which satisfy (\ref{eq102}), one may define $\theta_2$ by formula (\ref{eq101}), put $\theta _1 = \kappa\theta_2$ and define $\Psi$ by formula (\ref{eq4001}) (if only condition (\ref{eq411}) is satisfied), where $\xi _i = (1/f_i)\cdot(\partial/\partial x_i)$. 

Substitution of $\theta_2$ from (\ref{eq101}) and $\theta_1 = \kappa\theta_2$ produces from (\ref{eq403}) and (\ref{eq404}) 
the system of two equations with the unknown function $\Psi$. The integrability conditions for this new system of equations reads as
\begin{equation}\label{eq103}
2 (-\frac{(\partial_2 f_1)^2}{f_1^2}+\frac{\partial_2^2 f_1}{f_1}+\frac{\partial_1 f_2^2+f_2 (\partial_2 f_1\partial_1 k+k \partial_1\partial_2 f_1)}{f_2^2})=0
\end{equation}
and
\begin{eqnarray}\label{eq104}
&&\frac{2}{f_1^6 f_2^6} (-15 f_2^6 \partial_2 f_1(\partial_1 f_1)^3-f_1^4 f_2^2 \partial_2 f_1 (2 f_2^4 \partial_1 f_1\nonumber \\
&+&3 \partial_2 f_1\partial_2 f_2\partial_1 f_2+f_2 (2 \partial_2 k(\partial_2 f_1)^2-3 \partial_2^2 f_1\partial_1 f_2)) \nonumber\\
&-&f_1^6 (2 f_2^5 (\partial_2 k\partial_2 f_1+k \partial_2^2 f_1)+f_2^3 (3 \partial_2^2 k \partial_2^2 f_1+\partial_2 f_1\partial_2^3 k\nonumber \\
&+&3 \partial_2 k\partial_2^3 f_1+k \partial_2^4 f_1)+2 f_2^4 \partial_2 f_2\partial_1 f_2+15 (\partial_2 f_2)^3\partial_1 f_2\nonumber \\
&+&5 f_2 \partial_2 f_2 (3 \partial_2 k\partial_2 f_1\partial_2 f_2+3 k \partial_2 f_2\partial_2^2 f_1 \nonumber\\
&-&2 \partial_2^2 f_2\partial_1 f_2)-f_2^2 (12\partial_2 k\partial_2 f_2\partial_2^2 f_1+\partial_2 f_1 (6 \partial_2 f_2\partial_2^2 k\nonumber \\
&+&4 \partial_2 k\partial_2^2 f_2)+k (4\partial_2^2 f_1\partial_2^2 f_2+6 \partial_2 f_2\partial_2^3 f_1)-\partial_2^3 f_2\partial_1 f_2)) \nonumber\\
&+&f_1^5 f_2 (f_2 (\partial_2 f_1)^2 (15 \partial_2 k\partial_2 f_2-4 f_2 \partial_2^2 k)-11 f_2^2\partial_2 k\partial_2 f_1\partial_2^2 f_1 \nonumber\\
&-&3 k f_2^2 (\partial_2^2 f_1)^2+\partial_2 f_1 (8 f_2^4+18 (\partial_2 f_2)^2-5 f_2 \partial_2^2 f_2)\partial_1 f_2\nonumber \\
&+&f_2 (-21 \partial_2 f_2\partial_2^2 f_1+5 f_2 \partial_2^3 f_1) \partial_1 f_2+2 f_2^5\partial_1\partial_2 f_1)\nonumber \\
&+&f_1 f_2^5 \partial_1 f_1 (15 f_2 \partial_1 f_1\partial_1\partial_2 f_1+2 \partial_2 f_1 (9 \partial_1 f_1\partial_1 f_2\nonumber \\
&+&5 f_2 \partial_1^2 f_1))-f_1^2 f_2^4 (3 \partial_2 f_1\partial_1 f_1(\partial_1 f_2)^2+f_2 (-12 (\partial_2 f_1)^2 \partial_1 k\partial_1 f_1 \nonumber\\
&+&27 \partial_1 f_1\partial_1 f_2\partial_1\partial_2 f_1+\partial_2 f_1 (5 \partial_1 f_2\partial_1^2 f_1-6 \partial_1 f_1\partial_1^2  f_2)) \nonumber\\
&+&f_2^2 (4 \partial_1\partial_2 f_1\partial_1^2 f_1+6\partial_1 f_1\partial_1^2\partial_2 f_1+\partial_2 f_1\partial_1^3 f_1))\nonumber\\
&+&f_1^3 f_2^4 (6 (\partial_1 f_2)^2\partial_1\partial_2 f_1-2 (\partial_2 f_1)^2 (2\partial_1 k\partial_1 f_2+f_2 \partial_1^2 k)\nonumber \\
&+&6 f_2 \partial_1 f_2\partial_1^2\partial_2 f_1-\partial_2 f_1 (3 \partial_1 f_2\partial_1^2 f_2+f_2 (10\partial_1 k\partial_1\partial_2 f_1 +\partial_1^3 f_2))\nonumber\\ 
&+&f_2 (-6 k(\partial_1\partial_2 f_1)^2-3 \partial_1\partial_2 f_1\partial_1^2 f_2+f_2 \partial_1^3\partial_2 f_1))) = 0.
\end{eqnarray}

If equations (\ref{eq103}) and (\ref{eq104}) are satisfied, then the function $\Psi$ given by (\ref{eq4001}) is a solution to 
(\ref{eq403}) -- (\ref{eq404}) and  the system of quantities $(\xi_1, \xi_2, \theta_1, \theta _2, \Psi)$ satisfies 
integrability conditions of Section \ref{sec:lci}, so it corresponds to a unique (up to M\"obius transformation) surface $S$ 
for which the leaves of the Dupin foliation meet the lines of curvature at the angle $\alpha$ such that $\tan^3\alpha = \kappa$. 
This way, we proved the following.

\begin{thm}\label{thm2}
For any function $\kappa$ for which the system (\ref{eq102}), (\ref{eq103}), (\ref{eq104}) of partial differential equations 
(with unknown functions $f_1$ and $f_2$) has a solution, 
there exists a surface $S$ such that the Dupin lines on $S$ intersect lines of curvature at angle
$\alpha$ such that $\tan^3\alpha = \kappa$.\hfill\qed
\end{thm}

If our function $\kappa$ is constant, then equations (\ref{eq103}) and (\ref{eq104}) simplify, respectively, to
\begin{equation}\label{eq1003}
 \frac{2\partial_2^2 f_1}{f_1} + \frac{2k(-\partial_2 f_1 \cdot \partial_1 f_2+ f_2 \partial_1\partial_2 f_1)}{f_2^2} - \frac{2 (\partial_2 f_1)^2}{f_1^2} = 0
 \end{equation}
 and
 \begin{eqnarray}\label{eq1004}
 &-&\frac{1}{f_1^6 f_2^6} 2 (k f_1^6 (-2 f_2^4 \partial_2 f_1 \partial_2 f_2-15 \partial_2 f_1 (\partial_2 f_2)^3\nonumber \\
 &+&2 f_2^5 \partial_2^2 f_1+5 f_2 \partial_2 f_2 (3 \partial_2 f_2 \partial_2^2 f_1 \nonumber\\
 &+&2 \partial_2 f_1 \partial_2^2 f_2)-f_2^2 (4 \partial_2^2 f_1 \partial_2^2 f_2+6 \partial_2 f_2 \partial_2^3 f_1 \nonumber\\
 &+&\partial_2 f_1 \partial_2^3 f_2)+f_2^3 \partial_2^4 f_1)+15 f_2^6 \partial_2 f_1 (\partial_1 f_1)^3 \nonumber\\
 &+&f_1^4 f_2^2 \partial_2 f_1 (-3 k (\partial_2 f_1)^2 \partial_2 f_2+3 k f_2 \partial_2 f_1 \partial_2^2 f_1 \nonumber\\
 &+&2 f_2^4 \partial_1 f_1)+f_1^5 f_2 (8 k f_2^4 (\partial_2 f_1)^2+18 k (\partial_2 f_1)^2 (\partial_2 f_2)^2 \nonumber\\
 &-&k f_2 \partial_2 f_1 (21 \partial_2 f_2 \partial_2^2 f_1+5 \partial_2 f_1 \partial_2^2 f_2)\nonumber\\
 &+&k f_2^2 (3 (\partial_2^2 f_1)^2+5 \partial_2 f_1 \partial_2^3 f_1)-2 f_2^5 \partial_1\partial_2 f_1)\nonumber \\
 &+&f_1 f_2^5 \partial_1 f_1 (30 k (\partial_2 f_1)^2 \partial_1 f_1-15 f_2 \partial_1 f_1\partial_1\partial_2 f_1\nonumber \\
 &+&2 \partial_2 f_1 (6 \partial_1 f_1\partial_1 f_2-5 f_2 \partial_1^2 f_1)) \nonumber\\
 &+&f_1^2 f_2^4 (24 k^2 (\partial_2)^3 \partial_1 f_1+(\partial_2 f_1)^2 (30 k \partial_1 f_1 \partial_1 f_2\nonumber \\
 &-&8 k f_2\partial _1^2 f_1)+f_2 (4 f_2 \partial_1 \partial_2 f_1\partial_1^2 f_1 \nonumber\\
 &+&\partial_1 f_1 (-9 \partial_1 f_2\partial_1\partial_2 f_1+6 f_2 \partial _1^2\partial_2 f_1))\nonumber \\
 &+&\partial_2 f_1 (\partial_1 f_1 (9 (\partial_1 f_2)^2-6 f_2 (6 k \partial_1\partial_2 f_1+\partial_1^2 f_2) )\nonumber \\
 &+&f_2 (-3 \partial_1 f_2 \partial_1^2 f_1+f_2 \partial_1^3 f_1))) \nonumber\\
 &+&f_1^3 f_2^3 (8 k^3 (\partial_2 f_1)^4+20 k^2 (\partial_2 f_1)^3 \partial_1 f_2-8 k (\partial_2 f_1)^2 (-2 (\partial_1 f_2)^2\nonumber\\
 &+&f_2 (3 k \partial_1\partial_2 f_1+\partial_1^2 f_2))+\partial_2 f_1(4 (\partial_1 f_2)^3\nonumber \\
 &-&f_2 \partial_1 f_2 (22 k \partial_1\partial_2 f_1+5 \partial_1^2 f_2) \nonumber\\
 &+&f_2^2 (8 k \partial_1^2\partial_2 f_1+\partial_1^3 f_2))+f_2 (-4 (\partial_1 f_2)^2 \partial_1\partial_2 f_1 \nonumber\\
 &+&2 f_2 \partial_1 f_2 \partial_1^2\partial_2 f_1+f_2 (6 k (\partial_1\partial_2 f_1)^2 \nonumber\\
 &+&3 \partial_1\partial_2 f_1\partial_1^2 f_2-f_2 \partial_1^3\partial_2 f_1)))) = 0
\end{eqnarray}

Solving the system (\ref{eq102}), (\ref{eq1003}) and (\ref{eq1004}) seems to be still difficult but anyway we have the following. 

\begin{figure}[h!]
\begin{center}
\includegraphics[scale=0.5]{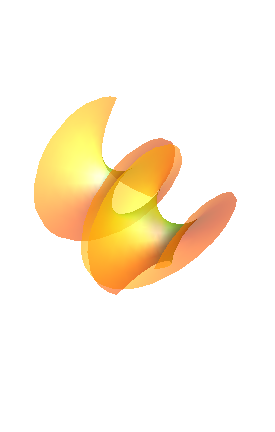}
\caption{A helcat\label{fig2}}
\end{center}
\end{figure}

\begin{ex}\label{ex2}{\rm The surfaces $S_\alpha$, $\alpha\in [0, \pi/2]$, see Figure \ref{fig2}, given by 
\begin{eqnarray*}
x_1&=& cos\alpha\cdot \sinh s\cdot \sin t + \sin\alpha\cdot\cosh s\cos t,\\
x_2 &=& - \cos\alpha\cdot\sinh s\cdot\cos t + \sin\alpha\cdot\cosh s\cdot\sin t,\\
x_3 &=& \sin\alpha\cdot s + \cos\alpha\cdot t,
\end{eqnarray*}
are called {\it helcats} (do not confuse them with an American movie!) and form a 1-parameter family of minimal surfaces connecting 
the helicoid $S_0$ with the catenoid $S_{\pi/2}$ (see, for example, \cite{op}). Since helicoids are the only ruled surfaces and catenoids the only surfaces of revolution among all the minimal surfaces, $S_\alpha$, $0 < \alpha < \pi/2$, is neither ruled nor a surface of revolution.
For these surfaces one has the following system of local conformal invariants:
\begin{equation*}
\omega_1= -\, {\frac {\sqrt {2}\cosh  s} {2 \sqrt {1+\sin\alpha }}}\cdot  \left( \cos\alpha
 \cdot ds - (\sin\alpha +1)\cdot dt \right) ,
\end{equation*}
 
 \begin{equation*}\omega_2= -{\frac {\sqrt {2}\cosh  s} {2 \sqrt {1+\sin\alpha }}}\cdot \left( \cos\alpha
 \cdot dt + (\sin\alpha +1)\cdot ds \right),
\end{equation*}
equivalently,
\begin{equation*}
\xi_1 = \frac{1}{\sqrt{2(1+\sin\alpha)}\cosh s}\cdot\Big(\cos\alpha\frac{\partial}{\partial s} 
- (1 + \sin\alpha)\frac{\partial}{\partial t}\big),
\end{equation*}

\begin{equation*}
\xi_2 = \frac{1}{\sqrt{2(1+\sin\alpha)}\cosh s}\cdot\Big(\cos\alpha\frac{\partial}{\partial t} 
+ (1 + \sin\alpha)\frac{\partial}{\partial s}\big),
\end{equation*}

\begin{equation*}\theta_1=\sqrt {2(1-\sin\alpha)}\cdot\sinh s, \quad
\theta_2= \sqrt {2(1+\sin\alpha)}\cdot\sinh s
\end{equation*}
and
\begin{equation*}
\Psi = \sin\alpha\cdot( 3\, \cosh^2 s  -2)
 \end{equation*}
Consequently,
\begin{equation*}\kappa = \theta_1/\theta_2 = {\frac {\cos\alpha }{1+\sin\alpha }}
\end{equation*} 
is constant on $S_\alpha$. In particular, $\kappa = 1$ on the helicoid and $\kappa = 0$ on the catenoid (and on all canal surfaces
as was mentioned before). Finally, recall that a surface $S$ is called {\it isothermic} if there exist on $S$ locally conformal 
parametrizations by curvature lines. This can be expressed by existence on $S$ of the nonzero function $u$ 
for which $[u\xi_1, u\xi_2] =0$. One can observe that the catenoid is the only surface among  $S_\alpha$'s which has this property. 

\begin{table}[h!]
\begin{center}
\begin{tabular}{|l|c|c|c|}
  \hline 
  $\alpha\ \ \ \ \ \ \backslash\ \ \ \ \ \  s$ & $0$ & $\pm 1$ & $\pm 2$ \\
  \hline
  0 & 1.5  & 1.5 & 1.5\\
  \hline
$\pi /100$&1.54&1.79&3.18\\
\hline
$\pi/7.384663\ldots$&2&5.12&31.7\\
\hline
$\pi/6$&2.07&5.84&37.96\\
\hline
$\pi/4$&2.17&7.44&52.34\\
\hline
$\pi/3$&2.12&8.44&62.33\\
\hline
$\pi/2.25$&1.82&8.6&66.32\\
\hline
$\pi/2.1$&1.68&8.29&64.62\\
\hline
$\pi/2.01$&1.53&1.79&61.68\\
\hline
\end{tabular} 
\end{center}
\caption{Values of $\Psi$ for cyclides osculating helcats.\label{tab1}}
\end{table}

One can ask also about the conformal type of the cyclide osculating helcats. Numerical experiments (performed with the use of Maple 14) show that
\begin{itemize}
\item for the helicoid, all the osculating cyclides are regular and have the same  invariant $\Psi$: $\Psi = 3/2$ at all the points,
\item for other helcats $S_\alpha$ with $\alpha$ positive and small enough, the osculating cyclides are regular along the axis $s = 0$ and become singular (first, just for one value of $s$ with one singularity, then, for larger values of $s$, with two singularities) as $s$ grows,
\item for helcats $S_\alpha$ with $\alpha$ large enough and reasonably smaller than $\pi/2$, all the osculating cyclides are singular and have two singularities.
\item as $\alpha$ approaches $\pi/2$, the osculating cyclides become again regular for small $s$ and still singular for $s$ large enough.
\end{itemize}
Approximate values of $\Psi$ for cyclides osculating helcats and for different values of $\alpha$ and $s$ (obviously, this value does not depend on the other parameter, $t$) are shown in Table \ref{tab1}.

Finally,  produced by Mathematica Figure \ref{fig4} shows the relative position of the osculating cyclide (in green) and the corresponding helcat $S_\alpha$ (in yellow), from the left to the right:
$\alpha = 0, \ \pi/4, \ \pi/2$, always at the point $s = t= 1$.

\begin{figure}[h!]
\begin{center}
\includegraphics[scale=0.4]{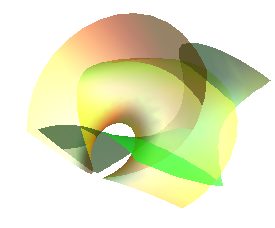}\quad 
\includegraphics[scale=0.4] {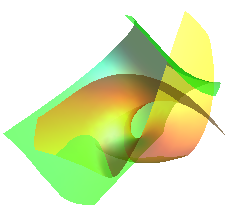}\quad 
\includegraphics[scale=0.4] {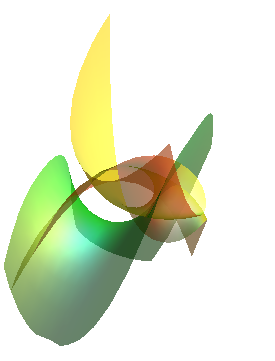}
\caption{Helcats and their osculating cyclides\label{fig4}}
\end{center}
\end{figure}
} \end{ex}


\begin{thebibliography}{WWW}

\bibitem[BLW]{blw} A. Bartoszek, R. Langevin, P. Walczak, {\it Special canal surfaces of \S3}, Bull. Braz. Math. Soc. {\bf 42} (2011), 301--320.

\bibitem[BW]{bw} A. Bartoszek, P. Walczak, {\it Foliations by surfaces of
a peculiar class}, Ann. Polon. Math., 94 (2008), 89 -- 95.

\bibitem[Br]{br} R. Bryant, {\it A duality theorem for Willmore surfaces},
J. Diff. Geom. {\bf 20} (1084), 23 -- 53.

\bibitem[CSW]{csw} G. Cairns, R. W. Sharpe and L. Webb. {\it Conformal
invariants for curves in three dimensional space forms}, Rocky Mountain J.
Math. {\bf 24} (1994), 933 -- 959.

\bibitem[Ca]{ca} A. Cayley, {\it On the cyclide}, Quart. J. Pure Appl. Math. {\bf 12} (1873), 148--163.

\bibitem[CG]{cg} R. Cipolla, P.  Giblin, {\it Visual motion of curves and surfaces},
Cambridge Univ. Press, Cambridge 2000.

\bibitem[Da1]{da1} G. Darboux, {\it Sur le contact des courbes et des surfaces}, Bull. sci. math. et astr., {\bf 4} (1880), 348 -- 384.

\bibitem[Da2]{da} G. Darboux, {\it Le\c cons sur la th\'eorie  g\'en\'erale des surfaces}, Guthier-Villars, Paris 1897.

\bibitem[Dr1]{dr1} J. Drach, {\it Sur les lignes d'osculation quadrique des surfaces}, C. R. Acad. Sci. Paris {\bf  224} (1947), 309 -- 312. 

\bibitem[Dr2]{dr2} J. Drach, {\it Détermination des lignes d'osculation quadrique (lignes de Darboux) sur les surfaces cubiques. Lignes asymptotiques de la surface de Bioche},
C. R. Acad. Sci. Paris {\bf 226} (1948), 1561-1564.

\bibitem[Du]{du} C. Dupin, {\it Applications de G\'eom\'etrie  et de M\'echanique}, Bachelier, Paris 1822.
 
\bibitem[Fi]{fi} A. Fialkov, {\it Conformal differential geometry of a
subspace}, Trans. Amer. Math. Soc. {\bf 56} (1944), 309 -- 433.

\bibitem[GLW]{glw} R. Garcia, R. Langevin, P. Walczak, {\it Dynamical behaviour of Darboux curves}, preprint, arXiv.0912.3749. (2009).

\bibitem[LW]{lw} R. Langevin, P. Walczak, {\it Conformal geometry of foliations}, Geom. Dedicata {\bf 132}, 135 -- 178.

\bibitem[Ma]{ma} J. C. Maxwel, {\it On the cyclide}, Quart. J. Pure Appl. Math., {\bf 9} (1868), 111--126,
l
\bibitem[Op]{op} J. Oprea, {\it Differential Geometry and its Applications}, Prentice Hall 1997.

\bibitem[Po]{po} I. R. Porteous, {\it Geometric differentiation. For the intelligence of curves and surfaces}, Cambridge Univ. Press, Cambridge 2001.

\bibitem[Sa1]{sa1} L. A.~Santal\'{o}, {\it Curvas extremales de la torsion total y curvas-D}, Publ. Inst. Mat. Univ. Nac. Litoral.  1941,  131--156.

\bibitem[Sa2]{sa2} L. A.~Santal\'{o}, {\it Curvas D sobre conos,} Select Works of L.A. Santal\'{o}, Springer Verlag 2009,  317-325. 

\bibitem[Tr]{tr} A. Tresse, {\it Sur les invariants diff\'erentiels d'une
surface par rapport aux transformations conformes de l'espace}, C.R. Acad.
Sci. Paris {\bf 114} (1892), 948 -- 950.
\end{thebibliography}
\end{document}